

\baselineskip=14pt
\parskip=10pt

\font\eightrm=cmr8 

\magnification=\magstephalf

\def\1{{\overline{1}}}
\def\2{{\overline{2}}}
\parindent=0pt
\overfullrule=0in

\def\frac#1#2{{#1 \over #2}}
\centerline
{\bf Two Motivated Concrete Proofs (much better than the usual one) that}
\centerline
{\bf the Square-Root of 2 is Irrational}
\bigskip
\centerline
{\it  Doron ZEILBERGER}
\bigskip
\qquad\qquad\qquad\qquad\qquad\qquad {\it Dedicated to Zvi Artstein, the Poincar\'e of the Early 21st Century}

{\bf Preface}

In his recent {\it masterpiece} [A], Zvika Artstein explains, so well, why most (in fact all!) students have trouble
understanding formal proofs in mathematics, and following formal logic. After all, Evolution did not prepare us
to do formal proofs, and be uptight about (alleged) {\it absolute certainty}. The straightjacket of {\it rigorous proof},
and the (to me abhorrent) {\it axiomatic method}, that is still the ruling dogma today, is a very questionable
gift bequeathed to us by the gang of Plato, Aristotle,  Euclid and their buddies.
Like abstinence as a form of birth-control, it was good for a while, but now there are much better ways to
control birth, and achieve mathematical knowledge!

And not only students have trouble! Even members of the illustrious {\it Coll\`ege de France} struggle with proofs by
contradiction (see [Z]). And indeed proofs by contradiction, and in general, 
abstract reasoning, are so counter-intuitive.

On page 53 of [A], Artstein rehashes the standard textbook proof that $\sqrt{2}$ is irrational, and much later,
in the very insightful (and depressing!) final chapter on math education, p. 400, there is a Socratic 
(imaginary, but very realistic) dialog between a student and a teacher, trying to come to grips with
the {\it pilpul} of that proof.

And indeed, {\it proofs by contradiction} are unsatisfying regardless of your mathematical sophistication or maturity,
and an {\it explicit}, constructive, proof is always preferable, whenever possible, and if it is not possible to
completely avoid contradiction, one should minimize it.

The panacea to the ills of mathematics education, so eloquently described by Artstein in that final chapter of [A], is
to do things {\it empirically}, using our beloved electronic servants (soon to be our masters) with several phases.
First an empirical ``proof'' for the first $1000$ cases, then also using {\it computer experiments}, to
design a proof that even Euclid would approve of, but much nicer and more concrete.

It is also a very good idea to teach our children how to program a computer, using, say, Python, or Maple, and that way they
would get used to logical thinking much more efficiently (and with fewer tears), then torturing them
with Euclidean-style axiomatic drivel. When you program, and debug, you get {\it instant feedback}, and
while evolution did not prepare us to reason logically, it did make us {\it adaptable}, and
get good at things that are fun to do.

\vfill\eject

{\bf The Empirical Proof}

Go to Maple, and type the following line:

{\tt min(seq(seq(abs(a**2-2*b**2),a=1..1000),b=1..1000));}

and {\it immediately} get the answer: $1$. So it turns out that the smallest that $|a^2-2b^2|$ can get, for
$1 \leq a,b \leq 1000$, is $1$, and it is quite plausible that it is {\it never} $0$.

Of course, this is not yet a formal proof. There exist simple diophantine equations whose smallest
solution is huge, but at least any student can grasp that $a^2-2b^2$, at least for $1 \leq a,b \leq 1000$, is never $0$.

{\bf First Complete Proof}

Next, one can look for those pairs $[a,b]$ that achieve the minimum

{\tt Z:=[]: for a from 1 to 1000 do for b from 1 to 1000 do if abs(a**2-2*b**2)=1 then Z:=[op(Z),[a,b]] fi: od:od: print(Z): } \quad \quad ,

and get, after one second, the output
$$
[[1, 1], [3, 2], [7, 5], [17, 12], [41, 29], [99, 70], [239, 169], [577, 408]] \quad ,
$$
and now we can look for a {\bf pattern}, and as Artstein emphasized, {\it looking for patterns}
was a great gift bestowed to us by Evolution, to the extent that we sometimes overdo it, and
find spurious patterns, by {\it overfitting}. 
But in this case anyone with IQ $\leq 100$
could not fail to see that
$$
1+1=2 \quad , \quad 3+2=5, \quad , \quad 7+5=12 \quad , \quad 17+12=29 \quad , \quad 41+29=70 \quad , \quad \dots \quad .
$$
In other words the $b$ of each pair $[a,b]$ that solves $a^2-2b^2= \pm 1$ is the sum of the
$a$ and $b$ of the previous pair. Also any one with IQ $\leq 105$ can see another pattern
$$
3-2=1 \quad, \quad 7-5=2  \quad, \quad 17-12=5  \quad, \quad 41-29=12  \quad, \quad 99-70=29 \quad  \dots
$$
Calling the $n$-th pair $[a_n,b_n]$ we have
$$
b_{n}=a_{n-1}+b_{n-1} \quad, \quad a_n-b_n=b_{n-1} \quad,
$$
So
$$
b_n=a_{n-1}+b_{n-1} \quad, \quad a_n=a_{n-1}+2b_{n-1} \quad ,
$$
and we found a way to generate (potentially!) {\it infinitely} many  solutions to the diophantine equation $a^2-2b^2= \pm 1$.
In other words, if $[a,b]$ is a solution then so is $[a+2b,a+b]$.

Having discovered this property, (experimentally!), one can prove it {\it rigorously}, 
by typing

{\tt expand((a+2*b)**2-2*(a+b)**2);} 

(or do it by hand, lazy bones!), and get that
$$
(a+2b)^2-2(a+b)^2=-(a^2-2b^2) \quad.
$$

But we can go backwards! Letting $c=a+2b, d=a+b$, we get 
$$
a=2d-c \quad, \quad b=c-d \quad .
$$
Note that $a+b=d<d+c$) and of course
$$
a^2-2b^2=-(c^2-2d^2) \quad,
$$
and now we can say that suppose that there are positive integers such that $a^2=2b^2$, then there would
be another pair $[c,d]$ with the same property but a smaller sum, and since $[2,1]$ is {\bf not} a solution,
qed. Of course, this is still, at the end, a proof by contradiction, but of a much gentler kind, and
it also taught us much more, namely how to construct infinitely many solutions of the so-called
Pell equation $a^2-2b^2=\pm 1$, and get terrific approximations $a/b$ for $\sqrt{2}$.

{\bf Remarks}: {\bf 1.} Ironically, at the end we got the algebraic version of what was (very possibly) the original, geometric, proof
of the Greeks, that is {\bf not} the standard proof, but a picture of the square of side-length $a$, and inside of it,
from two opposite corners, squares of side-length $b$. Their intersection is a middle square of side-length
$b-(a-b)=2b-a$ and the complement of their union are the two squares whose corners are the
other two corners, each of side-length $a-b$.

{\bf 2.} Of course, the above proof is modeled after {\it continued fractions}.

{\bf 3.} For an ultra-finitist like myself, neither the statement (taken literally), nor any of the proofs of the `fact' that
$\sqrt{2}$ is irrational, alias `never' rational, makes any sense, since it tacitly assumes that there is
an {\it infinite supply} of integers. But it is easy to make the statement meaningful.
The corrected statement is for {\it symbolic} positive integer $N$, $a^2-2b^2 \neq 0$ for all $1 \leq a,b \leq N$.
Now all the proofs make sense.

{\bf Second Complete Proof} (without computer)

Every rational number, alias, quotient of two positive integers, $a/b$, can be written
as
$$
\frac{a}{b} =q + \frac{r}{b} \quad,
$$
where $q$ is an integer and $r$ is less than $b$. If $r=0$ we have an integer, but
otherwise $b/r$ is larger than $1$ and we can keep doing it
$$
\frac{b}{r}= s+ \frac{t}{r}
$$
with $t$ less $r$, etc.

For example
$$
\frac{11}{4}=2+ \frac{3}{4} \quad ,
$$
$$
\frac{4}{3}= 1 + \frac{1}{3} \quad ,
$$
$$
\frac{3}{1}=3 \quad,
$$
sooner or later one gets an integer. You can let students do it for many examples to convince themselves that
for any ratio of integers this process eventually ends.

But now let's apply this process to $\sqrt{2}$
$$
\sqrt{2}= 1 + (\sqrt{2}-1) \quad
$$
$$
\frac{1}{\sqrt{2}-1}=\frac{\sqrt{2}+1}{(\sqrt{2}-1)(\sqrt{2}+1) }=\frac{\sqrt{2}+1}{1}=\sqrt{2}+1=2+ (\sqrt{2}-1) \quad.
$$
But now we are  back to $\sqrt{2}-1$ and it is clear intuitively (and rigorously) that we will never end, since
we are stuck in an infinite loop, hence $\sqrt{2}$ is {\bf not} rational.

{\bf Third Complete Proof} (with computer)

Teach them about the Maple command {\tt convert(., confrac)}, and try it out for many rational numbers and convince
them that you always get a finite answer.

Then type,

 {\tt Digits:=20:convert(evalf(1+sqrt(2)),confrac);}

 {\tt Digits:=100:convert(evalf(1+sqrt(2)),confrac);}

 {\tt Digits:=1000:convert(evalf(1+sqrt(2)),confrac);}

(The system variable {\tt Digits} in Maple indicates the precision of its floating-point calculations), and get
$$
[2, 2, 2, 2, 2, 2, 2, 2, 2, 2, 2, 2, \dots]
$$
as far as you can see, and anyone with IQ $\leq 90$ can guess that it goes for ever, i.e. that, conjecturally  for now,
$$
\sqrt{2}+1=2+1/(2+1/(2+1/2+ ...) \quad .
$$
In order to prove this conjecture, you start with the right side, and call it $y$.
By self-similarity, $y=2+1/y$, and solving this quadratic equation yields that indeed
$y=\sqrt{2}+1$, and once again we have a fully rigorous (but found by experimentation!) proof that
$\sqrt{2}+1$, and hence $\sqrt{2}$, is irrational, since its (simple) continued fraction is {\it infinite}, while
rational numbers' claim to fame is that they have a {\it finite} (simple) continued fraction.

{\bf Conclusion} 

The main problem with both contemporary mathematical research, and mathematical education, is
{\it inertia}. The great impact of computer-kind is only very slowly starting to kick-in, and
the  pernicious  pre-computer mentality, and the dogmatic, and outdated, insistence on
rigorous proof, and only accepting as mathematical knowledge fully (allegedly!) rigorously-proved
results, where there are so many true, but not yet proved, results that are {\it obviously} true,
since they have very convincing heuristic proofs. For example, the Riemann Hypothesis, $P \neq NP$,
Goldbach's conjecture, and 
the facts that $e+\pi$ and $\zeta(5)$ are irrational numbers, are all {\it certainly} true, at least
more certain than Fermat's Last Theorem, whose truth is not at all obvious on heuristic grounds, and
whose claimed proof, by Andrew Wiles, is very complicated, and has only been checked by a few experts.

Of course, sometimes rigorous proofs {\it are} nice, and {\it do} give additional insight, but these
should become optional! In fact, algorithms are much more fun, and 
pedagogically, as {\it mental gymnastics}, to train students in logical thinking, should be prefered to proofs.

[{\eightrm Speaking of algorithms, above I used the Euclidean algoritm $a=qb+r$, so Euclid was not all bad!}]

{\bf Disclaimer}: All the opinions in this article are my own, and while they are inspired by
Artstein's masterpiece, are probably not completely shared by him.

\bigskip
{\bf References}

[A] Zvi Artstein, ``{\it Mathematics and the Real World}: The Remarkable Role of {\bf Evolution} in the Making of {\bf Mathematics}'',
Prometheus Books, 2014.

[Z] Doron Zeilberger, {\it Opinion 70: A Case Study in Math (and Logic!) Abuse by a Social Scientist: Jon Elster's Critique of Backwards Induction},
Dec. 21, 2005, \hfill\break
{\tt http://www.math.rutgers.edu/{}\~{}zeilberg/Opinion70.html} \quad .

\bigskip
\hrule
\bigskip
Doron Zeilberger, Department of Mathematics, Rutgers University (New Brunswick), Hill Center-Busch Campus, 110 Frelinghuysen
Rd., Piscataway, NJ 08854-8019, USA. \hfill\break
{\tt zeilberg at math dot rutgers dot edu}   $\,$ , $\,$ {\tt http://www.math.rutgers.edu/\~{}zeilberg/} $\,$ .

\bigskip
\hrule
\bigskip
Oct. 7, 2014.
\end